\documentclass[12pt,draft]{extartic}

\usepackage{amsmath,amsthm,amssymb,calc}
\usepackage[dvips]{graphics, graphicx}

\usepackage[english]{babel}

\def\eps{\varepsilon}

\newcounter{num}[section]

\newcommand{\Th}{\refstepcounter{num}
{\bf Theorem \arabic{section}.\arabic{num} }}

\newcommand{\Lemma}{\refstepcounter{num}
{\bf Lemma \arabic{section}.\arabic{num} }}

\newcommand{\Note}{\refstepcounter{num}
{\it Note \arabic{section}.\arabic{num} }}

\newcommand{\Proof}{{\bf Proof. }}

\def\eps{\varepsilon}
\def\_phi{\varphi}

\def\la{\lambda}

\def\F{\widehat}
\def\L{\Lambda}
\def\m{\times}

\def\ov{\overline}

\def\C{{\mathbb C}}

\def\Z_N{{\mathbb Z}_N}
\def\Z{{\mathbb Z}}

\def\Span{{\rm Span\,}}

\def\Gr{{\mathbf G}}

\def\dim{{\rm dim }}

\def\l{\left}
\def\r{\right}

\def\supp{{\rm supp\,}}

\author{Shkredov I.D. and Sergey Yekhanin}
\title{ Sets with large additive energy and
symmetric sets
\footnote{ The first author was supported Pierre Deligne's grant
based on his 2004 Balzan prize, President's of Russian Federation
grant N MK--1959.2009.1, grant RFFI N 06-01-00383 and grant
Leading Scientific Schools No. 691.2008.1. The second author is
with Microsoft Research Silicon Valley Lab., e-mail:
yekhanin@microsoft.com}
}
\date{}
\begin{document}
\maketitle

\begin{center}
 Annotation.
\end{center}

{\it \small
     We show that for any set $A$ in a finite Abelian group $\Gr$ that has at least $c |A|^3$ solutions to $a_1 + a_2 = a_3 + a_4$, $a_i \in A$ there exist sets $A' \subseteq A$ and $\L\subseteq \Gr$, $\L = \{ \la_1, \dots, \la_t \}$, $t \ll c^{-1} \log |A|$ such that $A'$ is contained in $\left\{\sum_{j=1}^t \eps_j \la_j\ |\ \eps_j \in \{ 0,-1,1 \}\right\}$ and $A'$ has $\gg c |A|^3$ solutions to $a'_1 + a'_2 = a'_3 + a'_4$, $a'_i \in A'$. We also study so--called symmetric sets or, in other words, sets of large values of convolution.
}
\\
%\\
%\\

\section{Introduction}
\label{sec:introduction}

Let $\Gr$ be a finite Abelian group. For sets $A,B \subseteq \Gr$
let $E(A,B)$ denote their {\it additive energy}
$$
    E(A,B) := \left|\{ a_1+b_1 = a_2+b_2 ~:~ a_1,a_2 \in A,\, b_1,b_2 \in B \}\right| \,.
$$
We write $E(A)$ in place of $E(A,A).$ For a set $\L = \{ \la_1,
\dots, \la_t \}$ in $\Gr$ let $\Span (\L)$ denote the set
$\left\{\left.\sum_{j=1}^t \eps_j \la_j\ \right|\ \eps_j \in
\{0,-1,1\}\right\}.$ A set $A$ satisfying $E(A) \ge c |A|^3$ for
some constant $c,$ is called a set of {\it large additive energy}.
Sets of large additive energy are very important in additive
combinatorics~\cite{Tao_Vu_book}. In \cite{Sanders_Sh} T. Sanders
obtained the following result about such sets.

\smallskip

\Th {\it
    Let $\Gr$ be a finite Abelian group, $A \subseteq \Gr$ be a set, and $c\in (0,1]$.
    Suppose $E (A) \ge c |A|^3;$ then there exist sets $A_1 \subseteq A$ and $\L\subseteq \Gr$, such that $|\L| \ll c^{-1} \log |A|,$
    $A_1 \subseteq \Span (\L)$ and $|A_1| \ge 2^{-2} c^{1/2} |A|.$
} \label{t:c_0.5}

A slightly weaker version of the theorem above (with with
$c^{1/2+\eps}$ instead of $c^{1/2}$) was obtained
in~\cite{Sh_doubling} using so--called $(C,\beta)$--connected
sets. In~\cite{Sanders_Sh} Sanders also considered a stronger
restriction on the set $A,$ namely $|A+A| \le c^{-1} |A|,$ and
obtained an improvement of theorem~\ref{t:c_0.5} in this setting
(see theorem~\ref{t:Sanders_SD_A+B} below). He also found an
interesting generalization of theorem~\ref{t:c_0.5} for the case
of two different sets $A$ and $B$.

\smallskip

\Th {\it
    Let $\Gr$ be a finite Abelian group, $A,B \subseteq \Gr$ be two sets, and $c\in (0,1]$.
    Suppose $|A + B| \le  c^{-1} |A|;$ then there is a set $\L\subseteq \Gr$, $|\L| \ll c^{-1} \log |A|$ such that $B \subseteq \Span (\L)$.
} \label{t:Sanders_SD_A+B}

\medskip

Applications of the theorem above can be found
in~\cite{Sanders_asymmetric}. In the current paper we obtain an
extension of Theorem \ref{t:c_0.5} for the case of two different
sets $A$ and $B.$ We also obtain a refinement of the theorem in
the case $A=B.$

\smallskip

\Th {\it
    Let $\Gr$ be a finite Abelian group, $A,B \subseteq \Gr$ be two sets, and $c\in (0,1]$. Suppose $E (A,B) \ge c |A| |B|^2;$ then there exist sets $B_1 \subseteq B$ and $\L\subseteq \Gr$, $|\L| \ll c^{-1} \log |A|$ such that $B_1 \subseteq \Span (\L)$ and
    \begin{equation}\label{f:E(A,B_1)}
        E (A,B_1) \ge 2^{-5} E(A,B).
    \end{equation}
    In particular, $|B_1| \ge 2^{-3} c^{1/2} |B|$.

} \label{t:E(A,B)}

\Note {
    The result above yields an improvement of Theorem \ref{t:c_0.5}. Indeed, suppose that in the previous theorem we have $B=A.$ Let $A_1 = B_1.$
    Then $E(A,A_1) \ge 2^{-5} E(A)$. Using the Cauchy--Schwartz inequality, we get
    $$
        E(A_1) \ge 2^{-10} E(A) \,.
    $$
    Therefore $|A_1| \ge 2^{-4} c^{1/3} |A|$ and the exponent is sharp (see below).

} \label{note:c^{1/3}}

\smallskip

In what follows we give three proofs of theorem~\ref{t:E(A,B)}. In
section~\ref{sec:proof} we give the first (Fourier analytic)
proof. In section~\ref{sec:large_convolutions} we establish a
result on large values of convolution of two sets (theorem
\ref{t:large_convolutions}). We then give the second proof of
theorem~\ref{t:E(A,B)}. Our proof relies on an idea of
Sanders~\cite{Sanders_asymmetric}. We do not use the Fourier
transform and get slightly weaker bounds. Later we generalize
theorem~\ref{t:large_convolutions}, and rely on that
generalization to obtain the last proof of (a small refinement of)
theorem~\ref{t:E(A,B)}. Again, we do not use the Fourier method.

Our results concerning the structure of sets of large values of
convolution are of independent interest. Our results on sets with
large additive energy are considerably weaker than the
implications of the polynomial Freiman--Ruzsa
conjecture~\cite{Green_FR_conj}.

\medskip

We conclude with few comments regarding the notation used in this
paper. For a positive integer $n,$ we set $[n]=\{1,\ldots,n\}.$
All logarithms are base $2.$ Signs $\ll$ and $\gg$ are the usual
Vinogradov's symbols. Finally, with a slight abuse of notation we
use the same letter to denote a set $S\subseteq \Gr$ and its
characteristic function $S:\Gr\rightarrow \{0,1\}.$

\medskip

The authors are grateful to T.~Sanders for useful discussions.
%Mathias Beiglb\"{o}ck, Alexander Fish for useful discussions
%and S.V. Konyagin for a number of helpful advices and remarks.
%Also I acknowledge the Mathematical Sciences Research Institute for its hospitality
%and providing me with excellent working conditions.

%\input{proof}

\section{Proof of the main result}
\label{sec:proof}

%First of all we obtain
%%Let us obtain an analog of Theorem \ref{t:c_0.5} for two different sets.

Let $\Gr$ be a finite Abelian group, $N=|\Gr|.$ It is
well--known~\cite{Rudin_book} that the dual group $\F{\Gr}$ is
isomorphic to $\Gr.$ Let $f$ be a function from $\Gr$ to $\C.$  We
denote the Fourier transform of $f$ by~$\F{f},$
\begin{equation}\label{F:Fourier}
  \F{f}(\xi) =  \sum_{x \in \Gr} f(x) e( -\xi \cdot x) \,,
\end{equation}
where $e(x) = e^{2\pi i x}$. We rely on the following basic
identities
\begin{equation}\label{F_Par}
    \sum_{x\in \Gr} \left|f(x)\right|^2
        =
            \frac{1}{N} \sum_{\xi \in \F{\Gr}} \left|\widehat{f} (\xi)\right|^2 \,.
\end{equation}

\begin{equation}\label{svertka}
    \sum_{y\in \Gr} \left|\sum_{x\in \Gr} f(x) g(y-x) \right|^2
        = \frac{1}{N} \sum_{\xi \in \F{\Gr}} \left|\widehat{f} (\xi)\right|^2 \left|\widehat{g} (\xi)\right|^2 \,.
\end{equation}
If
$$
    (f*g) (x) := \sum_{y\in \Gr} f(y) g(x-y)
$$
 then
\begin{equation}\label{f:F_svertka}
    \F{f*g} = \F{f} \F{g} \quad \mbox{ and } \quad (\F{fg}) (x) = \frac{1}{N} (\F{f} * \F{g}) (x) \,.
\end{equation}
Using (\ref{svertka}), we can express additive energy in terms of
the Fourier transform
$$
    E(A,B) = \frac{1}{N} \sum_{\xi} \left|\F{A} (\xi)\right|^2 \left|\F{B} (\xi)\right|^2 \,.
$$

Our first proof of theorem~\ref{t:E(A,B)} relies on the following
lemma of T.~Sanders~\cite{Sanders_Sh}. (Similar results were
obtained by J.~Bourgain~\cite{BourgainA+A} and by the first
author~\cite{Sh_doubling}.) Recall that a set $\L = \{ \la_1,
\dots, \la_t \}$ in a finite Abelian group $\Gr$ is called {\it
dissociated} if any identity of the form $\sum_{j=1}^t \eps_j
\la_j = 0 $, where $\eps_j \in \{0,-1,1\}$ yields $\eps_j = 0$,
$j\in [t].$

\smallskip

\Lemma {\it
    Let $\Gr$ be a finite Abelian group, $Q\subseteq \Gr$ be a set, $l$ be a positive integer. There is a set $Q_1 \subseteq Q$ such that all dissociated subsets of $Q_1$ have size at most $l$ and
    for all $p\ge 2$
    %we have
    the following holds
    \begin{equation}\label{}
        \l( \frac{1}{N} \sum_\xi \left| \F{Q} (\xi) - \F{Q}_1 (\xi) \right|^p \r)^{1/p}
            \ll
                \sqrt{p/l} \cdot |Q| \,.
    \end{equation}
} \label{l:Bourgain_diss_sets}

{\bf Proof of Theorem \ref{t:E(A,B)}} Apply Lemma
\ref{l:Bourgain_diss_sets} to the set $B$ with parameters
$p=2+\log |A|$ and $l =  \eta^{-1} c^{-1} \log |A|$, where $\eta
\in (0,1]$ is an appropriate constant that we fix later. Write
$\eps (x) = B(x) - B_1(x)$, where $B_1 \subseteq B$ is such that
all dissociated subsets of $B_1$ have size at most $l$. We have
$$
    N \cdot E (A,B) = \sum_{\xi} \left|\F{A} (\xi)\right|^2 \left|\F{B} (\xi)\right|^2
        =
            \sum_{\xi} \left|\F{A} (\xi)\right|^2 \left|\F{B}_1 (\xi)\right|^2
                +
$$
$$
                +
                    \l(
                        \sum_{\xi} \left|\F{A} (\xi)\right|^2 \ov{\F{B}_1 (\xi)} \F{\eps} (\xi)
                            +
                        \sum_{\xi} \left|\F{A} (\xi)\right|^2 \F{B}_1 (\xi) \ov{\F{\eps} (\xi)}
                    \r)
                        +
                            \sum_{\xi} \left|\F{A} (\xi)\right|^2 \left|\F{\eps} (\xi)\right|^2
                                =
$$
$$
                                =
                                    \sigma_0 + \sigma_1 + \sigma_2 \,.
$$

By the H\"{o}lder inequality, identity~(\ref{F_Par}), and our
choice of parameters, we have
\begin{equation}\label{tmp:18.02.2010}
    \sigma_2
        \le
            \l( \sum_{\xi} |\F{\eps} (\xi)|^{2p} \r)^{1/p}
                \cdot
                    \l( \sum_{\xi} |\F{A} (\xi)|^{\frac{2p}{p-1}} \r)^{1-1/p}
                        \ll
                            \frac{p}{l} |B|^2 |A| |A|^{1/p} N
                                \le
                                    2^{-1} c |A| |B|^2 N \,.
\end{equation}
Hence either $\sigma_0$ or $\sigma_1$ is at least $2^{-2} c |A|
|B|^2 N$. In the first case we are done. In the second case an
application of the Cauchy--Schwartz inequality yields
$$
    2^{-6} N^2 E^2(A,B)
        \le
            N \cdot E(A,B_1) \cdot \sigma_2 \,.
$$
Combining the inequality above with~(\ref{tmp:18.02.2010}) we
get~(\ref{f:E(A,B_1)}). This completes the proof of
Theorem~\ref{t:E(A,B)}.

\medskip

For a set $Q \subseteq \Gr$ let $\dim (Q)$ denote the size of the
largest dissociated subset of $Q.$ Clearly, for any set $Q
\subseteq \Gr$ there is a dissociated set $\L \subseteq Q$ such
that $|\L| = \dim (Q)$ and $Q \subseteq \Span (\L).$ Thus, all
theorems above can be viewed as results concerning the dimension
of certain subsets of sets with large additive energy.

\medskip

By note~\ref{note:c^{1/3}}, theorem~\ref{t:E(A,B)} yields an
improvement of theorem~\ref{t:c_0.5}. Nevertheless the method
from~\cite{Sh_doubling} is surprisingly sharp. Indeed, the
argument there proceeds in two steps. Firstly, one finds a
$(C,\beta)$--connected subset of $A$ of size approximately
$c^{1/2} |A|$ (see~\cite{Sh_doubling} for appropriate
definitions). Secondly, one proves that any connected set belongs
to a span of a set of size $O(c^{-1} \log |A|)$. It is not hard to
verify that the bound used on the second step is sharp. The
argument used on the first step also cannot be improved. We are
grateful to T. Sanders for pointing us to the following example
(see also \cite{Sh_exp2}, theorem 4.1).

Let $\Gr = (\Z/2\Z)^n.$ For a linear subspace $H$ of $\Gr,$ let
$\supp(H)=\{i\in [n]\ |\ \exists x\in H, x_i\ne 0\}$ denote the
support of $H.$ Set $A = \bigcup_{i=1}^t H_i$ to be a union of $t$
linear subspaces $\{H_i\}_{i\in [t]}$ that have the same size
$h\gg t$ and disjoint supports. It is not hard to show
(see~\cite{Sh_exp2} for details) that $t h^3 \le E(A) \ll t h^3 =
(th)^3 / t^2$ and any connected subset of $A$ has cardinality
$O(h) = O ((th) /(t^2)^{1/2}).$

\medskip

Observe that the exponent of $c$ in note~\ref{note:c^{1/3}} is the
best possible. Indeed, set $\Gr = (\Z/2\Z)^n,$ and set $A=H
\bigsqcup \L,$ where $H$ is a linear subspace of size
approximately $c^{1/3} |\L|,$ and $\L$ is a dissociated set. Now
%$E(H) \gg E(A)$
$E(A) \gg c |A|^3$ and for every set $A_1 \subseteq A$ such that
$\dim (A_1) \ll c^{-1} \log |A|,$ $|A_1| \ll c^{1/3} |A|$
necessarily holds.

\section{Large values of convolution}
\label{sec:large_convolutions}

The following theorem bounds the dimension of symmetric
sets~\cite{Tao_Vu_book}, or in other words, sets of large values
of convolution.

\smallskip

\Th {\it
    Let $\Gr$ be a finite Abelian group, $A,B \subseteq \Gr$ be two sets. Let $\sigma \ge 1$ be a positive real number. Finally, let
    $$
        S = \{ x \in \Gr ~:~ (A * (-B)) (x) \ge \sigma \} \,.
    $$
    Then
    \begin{equation}\label{f:large_convolutions}
        \dim (S) \ll \max \{ |A|, |B| \} \cdot \sigma^{-1} \cdot \log ( \min \{ |A|, |B| \} ) \,.
    \end{equation}
} \label{t:large_convolutions}
\\
\Proof Assume $|B|\leq |A|.$ Let $\L$ be the largest dissociated
subset of $S,$ $|\L|=\dim(S).$ Consider a simple bipartite graph
$G = (V,E)$ with parts $A$ and $B$ and colors $\lambda\in \L$ on
edges. A vertex $a\in A$ is connected to a vertex $b\in B$ by an
edge colored $\lambda\in \L$ if and only if $a-b=\la.$ Note that
all edges incident to a certain vertex have different colors. Also
note that $|E| \ge \sigma |\L|.$

For an edge $e\in E,$ let $\rm{col}(e)\in \L$ denote its color.
Let $C=\{e_1,\ldots,e_k\}\subseteq E^k$ be an arbitrary $k$-long
cycle in $G.$ We have
\begin{equation}
\label{Eqn:ClosedSum} \sum_{i\in [k]}(-1)^i\rm{col}(e_i)=0.
\end{equation}

Let $e_i, i\in [k]$ be an arbitrary edge of $C.$ We say that $e_i$
is a {\it special} edge, if for all $j\in [k]$ such that $i\ne j$
we have $\rm{col}(e_i)\ne \rm{col}(e_j).$ We say that $C$ is a
{\it special} cycle, if one (or more) of its edges are special.
Observe that if $C$ is a special cycle; then~(\ref{Eqn:ClosedSum})
gives a non-trivial dependence between the elements of $\L.$ Thus
to prove theorem~\ref{t:large_convolutions} it suffices to
establish the following

\smallskip

\Lemma {\it
    Suppose in the setting above we have $|\L| > 16 |A| \sigma^{-1} \log |B|;$ then there is a special cycle of length at most $4 \log |B|$ in $G.$
} \label{l:edges_colouring}

\smallskip

Our proof of lemma~\ref{l:edges_colouring} relies on the following
lemma of Erd\"{o}s~\cite{Handbook_CO}[p. 74, lemma 7.1].

\smallskip

\Lemma {\it
    Let  $\Gamma = (V,E)$ be a finite simple graph, $d$ be a positive integer, and $|E| > (d-1) |V|$. Then $\Gamma$ has a subgraph of minimum degree at least $d$.
} \label{l:Erdos_graph_deg}

{\bf Proof of Lemma \ref{l:edges_colouring}} We apply the
Erd\"{o}s' lemma to the graph $G$ to obtain a sub-graph
$G^\prime.$ Note that the degree of every vertex of $G^\prime$ is
at least $d=2^{-2} \sigma |\L| |A|^{-1}$. By the assumption of the
lemma we have $d>4\log |B|.$ To find a special cycle in
$G^\prime,$ we pick an arbitrary node $v_0 \in G^\prime\cap A$ and
start carefully constructing a binary sub-tree of $G^\prime$
rooted at $v_0.$

We assign every node in our sub-tree (other $v_0$) a color, which
is the color of the edge that comes from its parent.  We gradually
extend the depth of our binary tree trying to keep the following
invariant satisfied: ''{\it For every node $v$ in the tree: The
color of $v$ is different from the  colors of all ancestors and
siblings of ancestors of $v.$}"

Below is the pseudo-code of our tree construction procedure. Here
$Tree$ denotes the set of nodes that are already in the tree
(initially $Tree=\{v_0\}$). Further, for any $v\in Tree,$ $F(v)$
denotes the set of colors that includes the color of $v$ as well
as the colors of all ancestors and siblings of ancestors of $v.$
We repeat the following procedure incrementing the value of $i$
starting with $i=0:$

\smallskip

$
\begin{array}{l}
1. \ \ \mbox{{\bf For} every node $v$ at depth $i$ {\bf Do}} \\
2. \ \ \ \ \ \ \ \mbox{\bf{Begin}} \\
3. \ \ \ \ \ \ \ \ \ \ \ \ \ \ \mbox{Pick $v_1$ and $v_2$ to be two children of $v$ such that} \\
4. \ \ \ \ \ \ \ \ \ \ \ \ \ \ \ \ \ \ \ \ \ \ \ \ \ \mbox{col$(\{v,v_1\})\not \in F(v)$ and col$(\{v,v_2\})\not \in F(v)$} \\
5. \ \ \ \ \ \ \ \ \ \ \ \ \ \ \ \ \ \ \ \ \ \ \ \ \ \mbox{(If no two such children exist {\bf Abort.})} \\
6. \ \ \ \ \ \ \ \ \ \ \ \ \ \ \mbox{{\bf If} ($v_1\in Tree$) or ($v_2\in Tree$) {\bf Then Abort.}} \\
7. \ \ \ \ \ \ \ \ \ \ \ \ \ \ \ \ \ \ \ \ \ \ \ \ \ \ \ \ \mbox{{\bf Else}} \\
8. \ \ \ \ \ \ \ \ \ \ \ \ \ \ \ \ \ \ \ \ \ \ \ \ \ \ \ \ \ \ \ \ \ \ \ \ \mbox{{\bf Begin}} \\
9. \ \ \ \ \ \ \ \ \ \ \ \ \ \ \ \ \ \ \ \ \ \ \ \ \ \ \ \ \ \ \ \ \ \ \ \ \ \ \ \ Tree:=Tree\cup\{v_1,v_2\} \\
10. \ \ \ \ \ \ \ \ \ \ \ \ \ \ \ \ \ \ \ \ \ \ \ \ \ \ \ \ \ \ \ \ \ \ \ \ \ \ \  F(v_1):=F(v)\cup\left\{\rm{col}\left(\{v,v_1\}\right),\rm{col}\left(\{v,v_2\}\right)\right\} \\
11.  \ \ \ \ \ \ \ \ \ \ \ \ \ \ \ \ \ \ \ \ \ \ \ \ \ \ \ \ \ \ \ \ \ \ \ \ \ \ \ F(v_2):=F(v)\cup\left\{\rm{col}\left(\{v,v_1\}\right),\rm{col}\left(\{v,v_2\}\right)\right\} \\
12. \ \ \ \ \ \ \ \ \ \ \ \ \ \ \ \ \ \ \ \ \ \ \ \ \ \ \ \ \ \ \ \ \ \ \ \ \mbox{{\bf End}} \\
13. \ \ \ \ \ \ \ \rm{\bf{End}} \\
\end{array}
$

\smallskip

The lower bound on $d$ that we have implies that while we
construct the first $2\log |B|$ levels of our tree we will always
be able to find two edges emanating from a node that have suitable
colors.  (In other words, no abort on line 5 of the pseudo-code
will occur while $i\leq 2\log |B|.$) Now observe that all odd
depth nodes in the tree we construct belong to the set $B.$
Therefore our tree construction algorithm will necessarily
discover some cycle $C$ and abort (at line 6 of the pseudo-code)
at some depth $i \leq 2\log |B|.$ We claim that $C$ is special
cycle. Indeed, let $v_*$ be the node of the smallest depth in $C.$
It not hard to check that both edges incident to $v_*$ in $C$ are
special. This concludes the proof of lemma and
theorem~\ref{t:large_convolutions}.

\medskip

\Note An appropriate version of Chang's theorem
(see~\cite{Sanders_asymmetric} or~\cite{Sh_Rudin}) implies a bound
for $\dim (S)$ that is weaker than~(\ref{f:large_convolutions}).
Specifically, it yields
$$
    \dim (S) \ll |A||B| \cdot \sigma^{-2} \cdot \log ( \min \{ |A|, |B| \} ) \,.
$$

\smallskip

\Note { Inequality (\ref{f:large_convolutions}) is the best
possible. To see this let $\Gr = (\Z/2\Z)^n.$ Let $B$ be a
subspace, and let $A=B\dotplus\L$, where $\L$ is a dissociated
set. Now $\sigma \sim |B|$ and $\dim (S) \sim |\L| + \dim (B).$
One can get a similar example with $E(B) = o (|B|^3),$ setting
$A=H\dotplus \L_1 \dotplus \L_2$ and $B=H \dotplus \L_1,$ where
$\L_1,\L_2$ are dissociated sets and $H$ is a subspace (note that
by construction sets $A$ and $B$ are connected). }
\label{note:sharp_B*B}

\bigskip

We now proceed to the second

{\bf Proof of Theorem \ref{t:E(A,B)}}. Let
$$
    S_j = \{ x\in \Gr ~:~ c 2^{j-2} |B| \le (A*B) (x) < c 2^{j-1} |B| \}\,,
                j\in [s]\,, \quad s \ll \log (1/c) \,.
$$
By assumption $E (A,B) \ge c |A| |B|^2$. Hence
$$
    \sum_{j=1}^s \sum_{x\in S_j} (A*B)^2 (x) \ge 2^{-1} c |A| |B|^2 \,.
$$
Put $
    c_j = \frac{1}{|A| |B|^2} \sum_{x\in S_j} (A*B)^2 (x)
$. Then
\begin{equation}\label{tmp:02.03.2010_1}
    2^{-1} c \le \sum_{j=1}^s c_j \le c
\end{equation}
and by definition of $S_j$, we have $c_j \le c 2^{j-1}$. Fix $j\in
[s]$ such that $c_j \ge (2s)^{-1} c > 0$. We have
\begin{equation}\label{f:02.03.2010_1}
    \sum_{x\in S_j} (A*B) (x) = \sum_{x} B(x) (S_j * (-A)) (x)
        \ge
            2^{-j+1} \frac{c_j}{c} |A| |B| \,.
\end{equation}
Let
$$
    %B^{(j)}_2
    B_1
    = \{ x \in B ~:~ (S_j * (-A)) (x) \ge 2^{-j} c_j c^{-1} |A| \} \,.
$$
By Theorem \ref{t:large_convolutions} the following holds
$$
    %\dim (B^{(j)}_2)
    \dim (B_1)
        \ll \max \{ |S_j|, |A| \} \cdot |A|^{-1} \frac{2^{j} c}{c_j} \cdot \log |A| \,.
$$
Since $(c 2^{j-2})^2 |B|^2 |S_j| \le c_j |A| |B|^2$ it follows
that $|S_j| \le 16 \cdot 2^{-2j} c_j c^{-2} |A|$. If $\max \{
|S_j|, |A| \} = |S_j|;$ then
$$
 %   \dim (B^{(j)}_2)
    \dim (B_1)
        \ll
            c^{-1} 2^{-j} \log |A|
                \ll
                    c^{-1} \log |A| \,.
$$
Now consider
%the situation when
the case $\max \{ |S_j|, |A| \} = |A|$.
%Put $B_1 = B^{(j)}_2$.
We have
$$
%    \dim (B^{(j)}_2)
    \dim (B_1)
        \ll
            2^{j} s \log |A|
                \ll
                    c^{-1} \log (c^{-1}) \log |A| \,.
$$
Since
$$
    \sum_{x\in S_j} (A*B_1) (x)  \ge 2^{-j} c_j c^{-1} |A| |B|
$$
it follows that
$$
    \sum_{x\in S_j} (A*B_1) (x) (A*B) (x) \ge 2^{-2} c_j |A| |B|^2 \,.
$$
Here the definition of $S_j$ was used. By the Cauchy--Schwartz
inequality and the definition of $c_j$, we obtain
$$
    E(A,B_1) \ge 2^{-4} s^{-1} E(A,B) \gg \log^{-1} (c^{-1}) \cdot E(A,B) \,.
$$
This completes the proof.

\smallskip

We now generalize theorem~\ref{t:large_convolutions} to the case
of more than two sets.

\smallskip

\Th {\it
    Let $\Gr$ be a finite Abelian group, $k \ge 2$ be a positive integer, $A_1, \dots, A_k \subseteq \Gr$, $|A_1| \le |A_2| \le \dots \le |A_k|$ be sets, and $\sigma \ge 1$ be a real number. Let
    $$
        S = \{ x \in \Gr ~:~ (A_1 * \dots * A_{k-2} * A_k * (-A_{k-1})) (x) \ge \sigma \} \,.
    $$
    Then
    \begin{equation}\label{f:large_convolutions_many}
        \dim (S) \ll |A_1| \dots |A_{k-2}| |A_k| \cdot \sigma^{-1} \cdot \log |A_{k-1}| \,.
    \end{equation}
} \label{t:large_convolutions_many} \Proof Let $\L \subseteq S$ be
the maximal dissociated subset. Consider a simple bipartite graph
$G = (V,E)$ with parts $A:=A_1 \m \dots \m A_{k-2} \m A_k$ and
$B:=A_{k-1}$ and colors $\lambda\in \L$ on edges. A vertex
$(a_1,\dots, a_{k-2}, a_k) \in A$ is connected to a vertex $b\in
B$ by an edge colored $\lambda\in \L$ if and only if $a_1 + \dots
a_{k-2} + a_k - b=\la.$ Note that $|E| \ge \sigma |\L|.$ Also note
that all edges incident to a certain vertex $a\in A$ have
different colors. Finally observe that for any vertex $b\in B,$
there exist edges of at least $\deg(b) / (|A_1| \dots |A_{k-2}|)$
different colors that are incident to $b.$ The latter observation
follows from the fact that
\begin{equation}\label{f:max_la+b}
    \max_{\la \in \L,\, b\in B} (A_1 * \dots * A_{k-2} * A_k) (\la+b)
        \le
            |A_1| \dots |A_{k-2}| \,.
\end{equation}

\smallskip

To proceed we need in a simple generalization of the Erd\"{o}s
lemma.

\smallskip

\Lemma {\it
    Let  $\Gamma = (V,E)$ be a finite simple bipartite graph with parts $V_1$ and $V_2.$
    Suppose $d_1,d_2$ are positive integers such that $|E| > (d_1-1) |V_1| + (d_2-1) |V_2|;$
    then $\Gamma$ has a bipartite subgraph with parts $V'_1 \subseteq V_1$, $V'_2 \subseteq V_2$ such that
    $$
        \min_{v'_1 \in V'_1} \deg (v'_1) \ge d_1 \,, \quad \mbox{and} \quad
            \min_{v'_2 \in V'_2} \deg (v'_2) \ge d_2 \,.
    $$
} \label{l:Erdos_graph_deg} \Proof Take any minimal bipartite
subgraph $\Gamma' = (V',E')$ of $\Gamma$ such that $|E(\Gamma')| >
(d_1-1) |V_1 (\Gamma')| + (d_2-1) |V_2 (\Gamma')|$, where $V_1
(\Gamma') \subseteq V_1$, $V_2 (\Gamma') \subseteq V_2$ are the
parts of $\Gamma'.$ It is easy to see that $\Gamma'$ has the
required properties. This completes the proof of the lemma.

\smallskip

To prove theorem~\ref{t:large_convolutions_many} we apply the
generalized Erd\"{o}s' lemma to $G,$ and obtain a bipartite
subgraph $G^\prime$ with parts $A'\subseteq A$, $B'\subseteq B$
such that for all $a^\prime\in A^\prime$ and $b^\prime\in
B^\prime,$ $\deg (a') \ge 2^{-2} \sigma |\L| / (|A_1| \dots
|A_{k-2}| |A_k|)$ and $\deg (b') \ge 2^{-2} \sigma |\L| /
|A_{k-1}|.$ Next we apply the (tree construction) argument from
the proof of theorem~\ref{t:large_convolutions} to the graph
$G^\prime.$ It is not hard to see that argument yields a
non-trivial dependency between the elements of $\L$ provided
$$
    \frac{\sigma |\L|}{|A_1| \dots |A_{k-2}| |A_k|} \gg \log |A_{k-1}|
$$
and
$$
    \frac{\sigma |\L|}{|A_{k-1}| |A_1| \dots |A_{k-2}|} \gg \log |A_{k-1}|
$$
This concludes the proof.

\medskip

We now give our third proof of theorem~\ref{t:E(A,B)}. In fact we
prove a slightly stronger result (see the
inequality~(\ref{f:B_1_third}) below).

\smallskip

{\bf Proof of Theorem \ref{t:E(A,B)}}. Without a loss of
generality assume $|A| \ge |B|.$ By assumption $E(A,B) \ge c |A|
|B|^2.$ It follows that
\begin{equation}\label{f:B_1_third}
    \sum_{x} (B*A*(-A)) (x) B_1 (x) \ge 2^{-1} c |A| |B|^2 \,,
\end{equation}
where $B_1 = \{ x\in B ~:~ (B*A*(-A)) (x) \ge 2^{-1} c |A| |B|
\}.$ Theorem~\ref{t:large_convolutions_many} yields $\dim (B_1)
\ll c^{-1} \log |A|.$ Combining the inequality~(\ref{f:B_1_third})
and the Cauchy--Schwartz inequality, we get $E(A,B_1) \gg 2^{-2} c
|A| |B|^2$ and the theorem follows.

\bigskip

If in theorem~\ref{t:large_convolutions_many} some extra
information on the additive energy of the sets~$A_j$ is available;
then the bound~(\ref{f:large_convolutions_many}) can be refined
for $k\geq 3$ (see~\cite{Sh_Rudin}). The example of
Note~\ref{note:sharp_B*B} shows that the corresponding estimates
in~\cite{Sh_Rudin} are sharp.

\medskip

\Note Let $k$ be a positive integer and $\L = \{ \la_1, \dots,
\la_t \} \subseteq \Gr$ be a set. We say that $\L$ belongs to the
family $\mathbf{\L} (k)$ if any identity of the form
$$
        \sum_{j=1}^t \eps_j \la_j = 0\,, \quad \quad \eps_j \in \{0,\pm 1\}\,, \quad \quad \sum_{j=1}^t |\eps_j| \le k\,,
$$
yields $\eps_j = 0$, $j\in [t]$. For $E\subseteq \Gr$ let $\dim_k
(E)$ denote the cardinality of the largest subset of $E$ that
belongs to the family $\mathbf{\L} (k).$ We remark that the
results above will still hold if one replaces $\dim(S)$ with
$\dim_k(S)$, say, for $k=O(\log |\Gr|)$. (For
Theorem~\ref{t:E(A,B)} see~\cite{Sh_exp1}).

\end{document}